# DISCUSSION OF "LEAST ANGLE REGRESSION" BY EFRON ET AL.

By Jean-Michel Loubes and Pascal Massart

*Université Paris-Sud*

The issue of model selection has drawn the attention of both applied and theoretical statisticians for a long time. Indeed, there has been an enormous range of contribution in model selection proposals, including work by Akaike (1973), Mallows (1973), Foster and George (1994), Birgé and Massart (2001a) and Abramovich, Benjamini, Donoho and Johnstone (2000). Over the last decade, modern computer-driven methods have been developed such as All Subsets, Forward Selection, Forward Stagewise or Lasso. Such methods are useful in the setting of the standard linear model, where we observe noisy data and wish to predict the response variable using only a few covariates, since they provide automatically linear models that fit the data. The procedure described in this paper is, on the one hand, numerically very efficient and, on the other hand, very general, since, with slight modifications, it enables us to recover the estimates given by the Lasso and Stagewise.

**1. Estimation procedure.** The "LARS" method is based on a recursive procedure selecting, at each step, the covariates having largest absolute correlation with the response $y$. In the case of an orthogonal design, the estimates can then be viewed as an $l^1$-penalized estimator. Consider the linear regression model where we observe $y$ with some random noise $\varepsilon$, with orthogonal design assumptions:

$$y = X\beta + \varepsilon.$$

Using the soft-thresholding form of the estimator, we can write it, equivalently, as the minimum of an ordinary least squares and an $l^1$ penalty over the coefficients of the regression. As a matter of fact, at step $k = 1, \ldots, m$, the estimators $\hat{\beta}^k = X^{-1}\hat{\mu}^k$ are given by

$$\hat{\mu}^k = \arg \min_{\mu \in \mathbb{R}^n} (\|Y - \mu\|_n^2 + 2\lambda_n^2(k)\|\mu\|_1).$$

There is a trade-off between the two terms, balanced by the smoothing decreasing sequence $\lambda_n^2(k)$. The more stress is laid on the penalty, the more parsimonious the representation will be. The choice of the $l^1$ penalty enables us to keep the largest coefficients, while the smallest ones shrink toward zero in a soft-thresholding scheme. This point of view is investigated in the Lasso algorithm as well as in studying the false discovery rate (FDR).

So, choosing these weights in an optimal way determines the form of the estimator as well as its asymptotic behavior. In the case of the algorithmic procedure, the suggested level is the $(k+1)$-order statistic:

$$\lambda_n^2(k) = |y|_{(k+1)}.$$

As a result, it seems possible to study the asymptotic behavior of the LARS estimates under some conditions on the coefficients of $\beta$. For instance, if there exists a roughness parameter $\rho \in [0,2]$, such that $\sum_{j=1}^m |\beta_j|^\rho \leq M$, metric entropy theory results lead to an upper bound for the mean square error $\|\hat{\beta} - \beta\|^2$. Here we refer to the results obtained in Loubes and van de Geer (2002). Consistency should be followed by the asymptotic distribution, as is done for the Lasso in Knight and Fu (2000).

The interest for such an investigation is double: first, it gives some insight into the properties of such estimators. Second, it suggests an approach for choosing the threshold $\lambda_n^2$ which can justify the empirical cross-validation method, developed later in the paper. Moreover, the asymptotic distributions of the estimators are needed for inference.

Other choices of penalty and loss functions can also be considered. First, for $\gamma \in (0,1]$, consider

$$J_\gamma(\beta) = \sum_{j=1}^m (X\beta)_j^\gamma.$$

If $\gamma < 1$, the penalty is not convex anymore, but there exist algorithms to solve the minimization problem. Constraints on the $l^\gamma$ norm of the coefficients are equivalent to lacunarity assumptions and may make estimation of sparse signals easier, which is often the case for high-dimensional data for instance.

Moreover, replacing the quadratic loss function with an $l^1$ loss gives rise to a robust estimator, the penalized absolute deviation of the form

$$\tilde{\mu}^k = \arg \min_{\mu \in \mathbb{R}^n} (\|Y - \mu\|_{n,1} + 2\lambda_n^2(k)\|\mu\|_1).$$

Hence, it is possible to get rid of the problem of variance estimation for the model with these estimates whose asymptotic behavior can be derived from Loubes and van de Geer (2002), in the regression framework.

Finally, a penalty over both the number of coefficients and the smoothness of the coefficients can be used to study, from a theoretical point of view, the asymptotics of the estimate. Such a penalty is analogous to complexity penalties studied in van de Geer (2001):

$$\mu^\star = \arg\min_{\mu \in \mathbb{R}^n, k \in [1,m]} (\|Y - \mu\|_n^2 + 2\lambda_n^2(k)\|\mu\|_1 + \text{pen}(k)).$$

**2. Mallows' $C_p$.** We now discuss the crucial issue of selecting the number $k$ of influential variables. To make this discussion clear, let us first assume the variance $\sigma^2$ of the regression errors to be known. Interestingly the penalized criterion which is proposed by the authors is exactly equivalent to Mallows' $C_p$ when the design is orthogonal (this is indeed the meaning of their Theorem 3). More precisely, using the same notation as in the paper, let us focus on the following situation which illustrates what happens in the orthogonal case where LARS is equivalent to the Lasso. One observes some random vector $y$ in $\mathbb{R}^n$, with expectation $\mu$ and covariance matrix $\sigma^2 I_n$. The variable selection problem that we want to solve here is to determine which components of $y$ are influential. According to Lemma 1, given $k$, the $k$th LARS estimate $\widehat{\mu}_k$ of $\mu$ can be explicitly computed as a soft-thresholding estimator. Indeed, considering the order statistics of the absolute values of the data denoted by

$$|y|_{(1)} \geq |y|_{(2)} \geq \cdots \geq |y|_{(n)}$$

and defining the soft threshold function $\eta(\cdot, t)$ with level $t \geq 0$ as

$$\eta(x, t) = x \mathbb{1}_{|x| > t}\left(1 - \frac{t}{|x|}\right),$$

one has

$$\widehat{\mu}_{k,i} = \eta(y_i, |y|_{(k+1)}).$$

To select $k$, the authors propose to minimize the $C_p$ criterion

(2.1) $$C_p(\widehat{\mu}_k) = \|y - \widehat{\mu}_k\|^2 - n\sigma^2 + 2k\sigma^2.$$

Our purpose is to analyze this proposal with the help of the results on penalized model selection criteria proved in Birgé and Massart (2001a, b). In these papers some oracle type inequalities are proved for selection procedures among some arbitrary collection of projection estimators on linear models when the regression errors are Gaussian. In particular one can apply them to the variable subset selection problem above, assuming the random vector $y$ to be Gaussian. If one decides to penalize in the same way the subsets of variables with the same cardinality, then the penalized criteria studied in Birgé and Massart (2001a, b) take the form

(2.2) $$C'_p(\widetilde{\mu}_k) = \|y - \widetilde{\mu}_k\|^2 - n\sigma^2 + \text{pen}(k),$$

where $\text{pen}(k)$ is some penalty to be chosen and $\widetilde{\mu}_k$ denotes the hard threshold estimator with components

$$\widetilde{\mu}_{k,i} = \eta'(y_i, |y|_{(k+1)}),$$

where

$$\eta'(x,t) = x \mathbb{1}_{|x|>t}.$$

The essence of the results proved by Birgé and Massart (2001a, b) in this case is the following. Their analysis covers penalties of the form

$$\text{pen}(k) = 2k\sigma^2 C\left(\log\left(\frac{n}{k}\right) + C'\right)$$

[note that the FDR penalty proposed in Abramovich, Benjamini, Donoho and Johnstone (2000) corresponds to the case $C = 1$]. It is proved in Birgé and Massart (2001a) that if the penalty $\text{pen}(k)$ is heavy enough (i.e., $C > 1$ and $C'$ is an adequate absolute constant), then the model selection procedure works in the sense that, up to a constant, the selected estimator $\widetilde{\mu}_{\widetilde{k}}$ performs as well as the best estimator among the collection $\{\widetilde{\mu}_k, 1 \leq k \leq n\}$ in terms of quadratic risk. On the contrary, it is proved in Birgé and Massart (2001b) that if $C < 1$, then at least asymptotically, whatever $C'$, the model selection does not work, in the sense that, even if $\mu = 0$, the procedure will systematically choose large values of $k$, leading to a suboptimal order for the quadratic risk of the selected estimator $\widetilde{\mu}_{\widetilde{k}}$. So, to summarize, some $2k\sigma^2 \log(n/k)$ term should be present in the penalty, in order to make the model selection criterion (2.2) work. In particular, the choice $\text{pen}(k) = 2k\sigma^2$ is not appropriate, which means that Mallows' $C_p$ does not work in this context. At first glance, these results seem to indicate that some problems could occur with the use of the Mallows' $C_p$ criterion (2.1). Fortunately, however, this is not at all the case because a very interesting phenomenon occurs, due to the soft-thresholding effect. As a matter of fact, if we compare the residual sums of squares of the soft threshold estimator $\widehat{\mu}_k$ and the hard threshold estimator $\widetilde{\mu}_k$, we easily get

$$\|y - \widehat{\mu}_k\|^2 - \|y - \widetilde{\mu}_k\|^2 = \sum_{i=1}^{n} |y|^2_{(k+1)} \mathbb{1}_{|y_i| > |y|_{(k+1)}} = k|y|^2_{(k+1)}$$

so that the "soft" $C_p$ criterion (2.1) can be interpreted as a "hard" criterion (2.2) with random penalty

(2.3) $$\text{pen}(k) = k|y|^2_{(k+1)} + 2k\sigma^2.$$

Of course this kind of penalty escapes *stricto sensu* to the analysis of Birgé and Massart (2001a, b) as described above since the penalty is not deterministic. However, it is quite easy to realize that, in this penalty, $|y|^2_{(k+1)}$ plays the role of the apparently "missing" logarithmic factor $2\sigma^2 \log(n/k)$. Indeed, let us consider the pure noise situation where $\mu = 0$ to keep the calculations as simple as possible. Then, if we consider the order statistics of a sample $U_1, \ldots, U_n$ of the uniform distribution on $[0,1]$

$$U_{(1)} \leq U_{(2)} \leq \cdots \leq U_{(n)},$$

taking care of the fact that these statistics are taken according to the usual increasing order while the order statistics on the data are taken according to the reverse order, $\sigma^{-2}|y|^2_{(k+1)}$ has the same distribution as

$$Q^{-1}(U_{(k+1)}),$$

where $Q$ denotes the tail function of the chi-square distribution with 1 degree of freedom. Now using the double approximations $Q^{-1}(u) \sim 2\log(|u|)$ as $u$ goes to 0 and $U_{(k+1)} \approx (k+1)/n$ (which at least means that, given $k$, $nU_{(k+1)}$ tends to $k+1$ almost surely as $n$ goes to infinity but can also be expressed with much more precise probability bounds) we derive that $|y|^2_{(k+1)} \approx 2\sigma^2 \log(n/(k+1))$. The conclusion is that it is possible to interpret the "soft" $C_p$ criterion (2.1) as a randomly penalized "hard" criterion (2.2). The random part of the penalty $k|y|^2_{(k+1)}$ cleverly plays the role of the unavoidable logarithmic term $2\sigma^2 k \log(n/k)$, allowing the hope that the usual $2k\sigma^2$ term will be heavy enough to make the selection procedure work as we believe it does. A very interesting feature of the penalty (2.3) is that its random part depends neither on the scale parameter $\sigma^2$ nor on the tail of the errors. This means that one could think to adapt the data-driven strategy proposed in Birgé and Massart (2001b) to choose the penalty without knowing the scale parameter to this context, even if the errors are not Gaussian. This would lead to the following heuristics. For large values of $k$, one can expect the quantity $-\|y - \widehat{\mu}_k\|^2$ to behave as an affine function of $k$ with slope $\alpha(n)\sigma^2$. If one is able to compute $\alpha(n)$, either theoretically or numerically (our guess is that it varies slowly with $n$ and that it is close to 1.5), then one can just estimate the slope (for instance by making a regression of $-\|y - \widehat{\mu}_k\|^2$ with respect to $k$ for large enough values of $k$) and plug the resulting estimate of $\sigma^2$ into (2.1). Of course, some more efforts would be required to complete this analysis and provide rigorous oracle inequalities in the spirit of those given in Birgé and Massart (2001a, b) or Abramovich, Benjamini, Donoho and Johnstone (2000) and also some simulations to check whether our proposal to estimate $\sigma^2$ is valid or not.

Our purpose here was just to mention some possible explorations starting from the present paper that we have found very stimulating. It seems to us that it solves practical questions of crucial interest and raises very interesting theoretical questions: consistency of LARS estimator; efficiency of Mallows' $C_p$ in this context; use of random penalties in model selection for more general frameworks.

CNRS and Laboratoire de Mathématiques
UMR 8628
Equipe de Probabilités, Statistique
 et Modélisation
Université Paris-Sud, Bât. 425
91405 Orday cedex
France
e-mail: Jean-Michel.Loubes@math.u-psud.fr

Laboratoire de Mathématiques
UMR 8628
Equipe de Probabilités, Statistique
 et Modélisation
Université Paris-Sud, Bât. 425
91405 Orday cedex
France
e-mail: Pascal.Massart@math.u-psud.fr